  \documentstyle[12pt]{article}
\newtheorem{ournothing}{}[subsection]

\newtheorem{ourlemma}[ournothing]{Lemma}
\newtheorem{ourtheorem}[ournothing]{Theorem}

\newtheorem{ourcorollary}[ournothing]{Corollary}
\newtheorem{ourproposition}[ournothing]{Proposition}
\newtheorem{ourconjecture}[ournothing]{Conjecture}

\newcommand{\numero}[1]{
\addtocounter{section}{1}
\begin{center}{\bf \thesection .\
#1\vspace{-.1in}}\end{center}
\setcounter{subsection}{0}
\setcounter{subsubsection}{0}
\setcounter{ournothing}{0}\indent}

\newcommand{\subnumero}[1]{
\addtocounter{subsection}{1}
\bigskip

\noindent{\bf \thesubsection \ }
{\sc #1}
\setcounter{subsubsection}{0}
\setcounter{ournothing}{0}
\newline
\indent}

\newenvironment{parag}{
\addtocounter{subsubsection}{1}
\begin{ournothing}\rm }{
\end{ournothing}}

\newcommand{\eop}{\hfill $/$\hspace*{-.1cm}$/$\hspace*{-.1cm}$/$\vspace{.1in}}

\newcommand{\newparag}[1]{
\end{parag}
\begin{parag}
\label{#1}
}

\newenvironment{lemma}{
\addtocounter{subsubsection}{1}
\begin{ourlemma}}{
\end{ourlemma}}

\newenvironment{corollary}{
\addtocounter{subsubsection}{1}
\begin{ourcorollary}}{
\end{ourcorollary}}

\newenvironment{theorem}{
\addtocounter{subsubsection}{1}
\begin{ourtheorem}}{
\end{ourtheorem}}

\newenvironment{proposition}{
\addtocounter{subsubsection}{1}
\begin{ourproposition}}{
\end{ourproposition}}

\newenvironment{remark}{
\addtocounter{subsubsection}{1}
\begin{ournothing}\hspace*{-.2cm}---Remark: \rm }{
\end{ournothing}}

\newtheorem{ulemma}{Lemma}
\newtheorem{utheorem}[ulemma]{Theorem}

\newcommand{\rr}{{\bf R}}

\newcommand{\Gg}{{\cal G}}

\newcommand{\Xx}{{\cal X}}

\begin{document}

\section*{On the Breen-Baez-Dolan stabilization hypothesis for Tamsamani's
weak $n$-categories}

\noindent
Carlos Simpson\newline
CNRS, UMR 5580, Universit\'e Toulouse 3, 31062 Toulouse CEDEX, France.

\bigskip

In \cite{BaezDolan} Baez and Dolan established their {\em stabilization
hypothesis} as one of a list of the key properties that a good theory of higher
categories should have. It is the analogue for $n$-categories of the
well-known stabilization theorems in homotopy theory.

To explain the statement,
recall that Baez-Dolan introduce the notion of {\em $k$-uply monoidal
$n$-category}  which is an $n+k$-category having only one $i$-morphism for all
$i< k$. This includes the notions previously defined and examined by many
authors,
of monoidal (resp. braided monoidal, symmetric monoidal) category (resp.
$2$-category) and so forth, as is explained in \cite{BaezDolan} \cite{cat}.
See the bibliographies of those preprints as well as that of the
the recent preprint \cite{BreenRecent} for many references concerning these
types of objects.  In the case where the $n$-category in question is an
$n$-groupoid, this notion is---except for truncation at $n$---the same thing
as the notion of $k$-fold iterated loop space, or
``$E_k$-space'' which appears in Dunn \cite{Dunn}
(see also some anterior references from there).
The fully stabilized notion of $k$-uply monoidal
$n$-categories for $k\gg n$ is what Grothendieck calls {\em Picard
$n$-categories}
in \cite{PursuingStacks}.

The {\em stabilization
hypothesis} \cite{BaezDolan} states that for $n+2 \leq  k\leq k' $, the $k$-uply
monoidal $n$-categories are the same thing as the $k'$-uply monoidal
$n$-categories.

This statement first appeared in a preliminary way in Breen
\cite{BreenAsterisque}; also there is some related correspondence
between Breen and Grothendieck in \cite{PursuingStacks}.

We will consider this hypothesis for Tamsamani's theory of
(``weak'') $n$-categories \cite{Tamsamani}, and show one of the main parts of
the statement, namely that a $k$-uply monoidal $n$-category can be ``delooped''
to a $k+1$-uply monoidal $n$-category, when $k\geq n+2$.

Before giving the precise statement, we make a
change of indexing. A {\em $k$-connected $n$-category} is an $n$-category
which has up to equivalence only one $i$-morphism for each $i\leq n$.
More precisely this means that the {\em truncation} $\tau _{\leq k}(A)$ is
trivial, equivalent to $\ast$. Note that a
$k-1$-connected $n+k$-category is equivalent to a
$k$-uply monoidal $n$-category (see \ref{prooflemma2} below).

We prove the following theorem (Corollary \ref{thm1}) for the theory
of \cite{Tamsamani}:

\begin{utheorem}
\label{thmintro}
If $A$ is a $k$-connected weak $n$-category and if $2k\geq n$
then there is a $k+1$-connected weak $n+1$-category $Y$ (with an object we
denote
by $y\in Y_0$) together with an equivalence $A\cong Hom _Y(y,y)$.
\end{utheorem}

\medskip

Translated back into the notation of \cite{BaezDolan} this says that
if $A$ is a $k$-uply monoidal  $n$-category and if $k\geq n+2$
then there exists a ``delooping'' $Y$ of $A$ which is a $k+1$-uply monoidal
$n$-category.

The statement of this theorem is the main content of the ``Stabilization
hypothesis'', but this result still leaves much further work to be done, for
example one would like to show that the stabilization construction $A\mapsto Y$
induces an equivalence of the higher categories which parametrize these objects
(and in fact one of the problems here is to pin down the right definition of
these parametrizing categories).

The technique we use for proving Theorem \ref{thmintro} is to remark that
one can reason with ``dimensions of cells'' for $n$-categories, in exactly the
same way as for topological spaces. The eventual non-invertibility of the
$i$-morphisms up to equivalence doesn't interfere. With this line of reasoning
available, the same arguments as in the classical topological case work.

Heavy reference will be made specially to notations from \cite{limits} and
to the closed model structure of \cite{nCAT}. For the reader's convenience I
have recopied (almost verbatim) several sections from \cite{limits} in the first
section of this note. There are some changes with respect to \cite{limits}
and notably the correction of an erratum in the definition of $\Upsilon$.

In \S 3, we will give some examples. These are mostly from
Baez-Dolan \cite{BaezDolan} and \cite{cat}, including: monoidal or symmetric
monoidal objects in $n-Cat$ giving rise to $n+1$ or $n+k$-categories; the
Whitehead operation which is useful for obtaining $k-1$-connected
$n$-categories;
the ``generalized center'' \cite{cat}; suspension; and the free $k$-uply
monoidal
$n$-category on one generator. In the last subsection we treat a possibly new
construction called ``cohomological twisting''. It is a generalization of the
classical construction in gauge theory whereby a gauge group $\Gg$ is
replaced by
its twist $Ad(P)$ with respect to a principal bundle $P$.

The present note was occasionned by the preprint ``Categorification'' by
Baez and
Dolan \cite{cat}.

\bigskip

\numero{Preliminary remarks}

For the reader's convenience, this first section consists of some  explanatory
remarks about Tamsamani's weak $n$-categories and related notions. These remarks
are mostly CUT-AND-PASTED out of \cite{limits}. The reader is strongly
advised to
refer to \cite{Tamsamani}, \cite{nCAT}, \cite{limits} and \cite{descente} for
further explanation. The reader already familiar with those references may skip
directly to the second section. Note however that we correct an {\em erratum}
in the definition of $\Upsilon$ in \cite{limits}, and add a small number of
remarks to the presentation.

\pagebreak[4]

\subnumero{$n$-categories}

\begin{parag}
\label{catnerve}
We begin by recalling the correspondence between categories and their nerves.
Let $\Delta$ denote the simplicial category whose objects are finite ordered
sets $p= \{ 0,\ldots , {\rm p}\}$ and morphisms are order-preserving maps.
If $C$ is a category then its nerve is the simplicial set (i.e.
a functor $A:\Delta ^o\rightarrow Sets$) defined by setting $A_p$ equal to the
set of composable $p$-uples of arrows in $C$.  This satisfies the property that
the ``Segal maps'' (cf the discussion of Segal's delooping machine \cite{Segal}
in \cite{Adams} for the origin of this terminology)
$$
A_p \rightarrow A_1 \times _{A_0} \ldots \times _{A_0} A_1
$$
are isomorphisms. To be precise this  map is given by the $p$-uple
of face maps $1\rightarrow p$ which take $0$ to ${\rm i}$ and $1$ to ${\rm i}+1$
for ${\rm i}=0,\ldots , {\rm p}-1$. Conversely, given a simplicial set $A$ such
that the Segal maps are isomorphisms we obtain a category $C$ by taking
$$
Ob (C) := A_0
$$
and
$$
Hom _C(x,y):= A_1(x,y)
$$
with the latter defined as the inverse image of $(x,y)$ under the map (given
by the pair of face maps) $A_1\rightarrow A_0 \times A_0$.  The condition on
the Segal maps implies that (with a similar notation)
$$
A_2(x,y,z)\stackrel{\cong}{\rightarrow} A_1(x,y)\times A_1(y,z)
$$
and the third face map $A_2(x,y,z)\rightarrow A_1(x,z)$ thus provides the
composition of morphisms for $C$. By looking at $A_3(x,y,z,w)$ one sees that the
composition is associative and the degeneracy maps in the simplicial set
provide the identity elements.
\end{parag}

\begin{parag}
\label{ncatsdef0}
Tamsamani's notion of weak $n$-category \cite{Tamsamani}
is a generalization of the above point of view on categories. We present
his definition in a highly recursive way, using the notion of
$n-1$-category in the definition of $n$-category. The original approach of
\cite{Tamsamani} was more direct. This definition is based on Segal's
delooping machine \cite{Segal} \cite{Adams}. It is also related to a definition
of Dunn \cite{Dunn} as is explained in \cite{descente}.

The related notion of ``Segal category'', used by the author in
\cite{effective}, was actually first mentionned by Dwyer, Kan and Smith
(\cite{DKS}, 1986) who prove the equivalence of the homotopy theories of
simplicially enriched categories and Segal categories.
\end{parag}

\begin{parag}
\label{notstrict}
Note that Tamsamani uses the terminology {\em $n$-nerve} for what we will call
``$n$-category'' since he needed to distinguish this from the notion of
strict $n$-category.  In the present paper we will (almost) never speak of
strict
$n$-categories and our terminology ``$n$-category'' means weak $n$-category or
$n$-nerve in the sense of \cite{Tamsamani}.
\end{parag}

\begin{parag}
\label{ncatsdef1}
 An {\em $n$-category} according to \cite{Tamsamani} is a functor $A$
from $\Delta ^o$ to the category of $n-1$-categories denoted
$$
p\mapsto A_{p/}
$$
such that $0$ is mapped to a set
\footnote{
Recursively an $n$-category which is a set is  a constant functor
where the $A_{p/}$ are all the same set---considered as $n-1$-categories.}
$A_0$
and such that the {\em Segal maps}
$$
A_{p/} \rightarrow A_{1/} \times _{A_0} \ldots \times _{A_0}A_{1/}
$$
are equivalences of $n-1$-categories (cf \ref{defequiv1} below).
\end{parag}

\begin{parag}
\label{multisimplicial}
The {\em category of $n$-categories} \cite{Tamsamani} denoted $n-Cat$ is
just the
category whose objects are as above and whose morphisms are the morphisms
strictly preserving the structure.
It is a subcategory of $Hom (\Delta ^o, (n-1)-Cat)$. Working this out
inductively we find in the end that $n-Cat$ is a subcategory of
$Hom ((\Delta ^n)^o, Sets)$, in other words an $n$-category is a
certain kind of multisimplicial set.
The multisimplicial set is denoted
$$
(p_1,\ldots , p_n )\mapsto A_{p_1,\ldots
p_n}
$$
and the $(n-1)$-category $A_{p/}$  itself considered as a multisimplicial
set has the expression
$$
A_{p/} = \left( (q_1,\ldots , q_{n-1})\mapsto A_{p,q_1,\ldots ,
q_{n-1}}\right) .
$$
\end{parag}

\begin{parag}
\label{theta}
The condition that $A_0$ be a set yields by induction the condition that
if $p_i=0$ then the functor $A_{p_1,\ldots , p_n}$ is independent of the
$p_{i+1}, \ldots , p_n$ cf
\cite{Tamsamani}. We call this the {\em constancy condition}. In
\cite{nCAT} we introduce the category $\Theta ^n$ which is the quotient of
$\Delta ^n$ defined by the condition that functors $(\Theta ^n)^o\rightarrow
Sets$ correspond to functors on $\Delta ^n$ having the above constancy
property.
Now $n-Cat$ is a subcategory of the category of presheaves of sets on $\Theta
^n$.  \end{parag}

\begin{parag}
\label{theta2}
To be more precise,  $\Theta ^n$ is defined to
be the quotient of the cartesian product $\Delta ^n$ obtained by identifying
all of the objects $(M, 0, M')$ for fixed $M = (m_1,\ldots , m_k)$ and variable
$M'= (m'_1, \ldots , m'_{n-k-1})$.  The object of $\Theta ^n$ corresponding
to the class of $(M,0,M')$ with all $m_i >0$ will be denoted
$M=(m_1,\ldots , m_k)$ with $k\leq n$
being the {\em length} denoted also $k=|M|$.
There is a unique object of length zero denoted $0$ by convention.
Two morphisms from $M$ to $M'$ in $\Delta ^n$ are identified if they both factor
through something of the form $(u_1,\ldots , u_i, 0, u_{i+2}, \ldots , u_n)$
and if their first $i$ components are the same.
\end{parag}

\begin{parag}
\label{arrows}
Let $1^i:= (1,\ldots , 1)\in \Theta ^n$. Recall from \cite{Tamsamani} that an
{\em $i$-arrow} in an $n$-category $A$ means an element of the set $A_{1^i}$.
There are notions of ``source'' and ``target'' of an $i$-arrow which are
$i-1$-arrows.
\end{parag}

\begin{parag}
Before discussing the notion of equivalence which enters into the above
definition of $n$-category we take note of the relationship with
\ref{catnerve}.
If $A$ is an $n$-category then its {\em set of objects} is the set $A_0$.  The
face maps give a morphism from $n-1$-categories to sets
$$
A_{p/}\rightarrow A_0 \times \ldots \times A_0
$$
and  we denote by $A_{p/}(x_0,\ldots , x_p)$ the $n-1$-category
inverse image of $(x_0,\ldots , x_p)$ under this map. For two objects $x,y\in
A_0$ the $n-1$-category $A_{1/}(x,y)$ is the {\em $n-1$-category of morphisms
from $x$ to $y$}.  This is the essential part of the structure which
corresponds, in the case of categories, to the $Hom$ sets.  One could adopt
the notation
$$
Hom _A(x,y):= A_{1/}(x,y).
$$
The condition that the Segal maps are equivalences of
$n-1$-categories says that the  $A_{p/}(x_0,\ldots , x_p)$ are determined up to
equivalence by the $A_{1/}(x,y)$. The role of the higher $A_{p/}(x_0,\ldots ,
x_p)$ is to  provide the composition (in the case $p=2$) and to keep track of
the higher homotopies of associativity $(p\geq 3$). Contrary to the case of
$1$-categories, here we need to go beyond $p=3$.
\end{parag}

\begin{parag}
\label{defequiv1}
In order for the recursive definition of $n$-category given in \ref{ncatsdef1}
to make sense, we need to know what an {\em equivalence} of $n$-categories is.
For this we generalize the usual notion for categories: an equivalence of
categories is a morphism which is (1) fully faithful and (2) essentially
surjective. We would like to define what it means for a functor between
$n$-categories $f:A\rightarrow B$ to be an equivalence. The generalized
{\em fully faithful condition} is immediate: we require that for any objects
$x,y\in A_0$ the morphism
$$
f: A_{1/} (x,y) \rightarrow B_{1/}(f(x), f(y))
$$
be an equivalence of $n-1$-categories (and we know what that
means by induction).
\end{parag}

\begin{parag}
\label{essentialsurjectivity}
The remaining question is how to define the notion
of essential surjectivity.  Tamsamani does this by defining a truncation
operation $T$ from $n$-categories to $n-1$-categories (a generalization of the
truncation of topological spaces used in the Postnikov tower).  Applying this
$n$ times to an $n$-category $A$ he obtains a set $T^nA$ which can also be
denoted $\tau _{\leq 0}A$.  This set is the set of ``objects of $A$ up to
equivalence'' where equivalence of objects is thought of in the
$n$-categorical sense.  We say that $f:A\rightarrow B$ is {\em essentially
surjective} if the induced map
$$
\tau _{\leq 0} (f) : \tau _{\leq 0} A \rightarrow \tau _{\leq 0} B
$$
is a surjection of sets. One has in fact that if $f$ is an equivalence according
to the above definition then $\tau _{\leq 0} f$ is an isomorphism.
\end{parag}

\begin{parag}
\label{anotherapproach}
Another way to approach the definition of
$\tau _{\leq 0}A$ is by induction in the following way. Suppose we know what
$\tau _{\leq 0}$ means for $n-1$-categories. Then for an $n$-category $A$ the
simplicial set $p\mapsto \tau _{\leq  0} (A_{p/})$ satisfies the condition that
the Segal maps are isomorphisms, so it is the nerve of a $1$-category. This
category may be denoted $\tau _{\leq 1} A$.  We then define $\tau _{\leq 0}A$
to be the set of isomorphism  classes of objects in the $1$-category
$\tau _{\leq 1} A$.
\end{parag}

\begin{parag}
\label{highertau}
We define similarly the truncation $\tau _{\leq k}(A)$ of an $n$-category $A$,
which is a $k$-category. This may be defined inductively (starting with the
knowledge of $\tau _{\leq 0}$) by the formula
$$
\tau _{\leq k}(A) _{p/} := \tau _{\leq k-1}(A_{p/}).
$$
In terms of Tamsamani's notation for truncation \cite{Tamsamani},
$$
\tau _{\leq k}(A)= T^{n-k}(A).
$$
\end{parag}

\bigskip

\subnumero{The closed model structure}

An $n$-category is a presheaf of sets on $\Theta ^n$ (\ref{theta}) satisfying
certain conditions as described above. Unfortunately $n-Cat$ considered as a
subcategory of the category of presheaves, is not closed under pushout or fiber
product. This remark is the starting point for \cite{nCAT}.  There, one
considers the full category of presheaves of sets on $\Theta ^n$ (these
presheaves
are called {\em $n$-precats}) and \cite{nCAT} provides a closed model structure
(cf \cite{Quillen} \cite{QuillenAnnals} \cite{Jardine})
on the category $nPC$ of $n$-precats, corresponding to
the homotopy theory of $n$-categories. In this section we briefly recall how
this works.

\begin{parag}
\label{precat}
An {\em $n$-precat} is defined to be a presheaf on the category $\Theta ^n$.
This corresponds to an $n$-simplicial set $(\Delta ^n)^o\rightarrow Sets$
which satisfies the constancy condition (cf \ref{theta}). The category $nPC$
of $n$-precats (with morphisms being the morphisms of presheaves) is to be given
a closed model structure.
\end{parag}

\begin{parag}
\label{precat2}
As in \ref{ncatsdef1}, an $n$-precat may also be viewed as a simplicial object
in the category $(n-1)PC$ satisfying the constancy condition that the zero$th$
term is a set.
\end{parag}

\begin{parag}
\label{prelims}
Note for a start that $nPC$ is closed under arbitrary products and coproducts,
also it admits an internal
$\underline{Hom}(A,B)$. These statements  come simply from the fact that $nPC$
is a category of presheaves over a category $\Theta ^n$.

We denote the coproduct or pushout of $A\rightarrow B$ and $A\rightarrow C$
by $B\cup ^AC$. We denote fiber products by the usual notation.
\end{parag}

\begin{parag}
\label{cofibs}
{\em Cofibrations:}
A morphism $A\rightarrow B$ of $n$-precats is a {\em cofibration} if the
morphisms $A_M \rightarrow B_M$ are injective whenever $M\in \Theta ^n$ is
an object of non-maximal length, i.e. $M= (m_1,\ldots , m_k, 0,\ldots , 0)$
for $k< n$.  The case of sets ($n=0$) shows that we can't require
injectivity at the top level $n$, nor do we need to.

We often use the notation $A\hookrightarrow B$ for a cofibration, not meaning
to imply injectivity at the top level.
\end{parag}

\begin{parag}
\label{opcat}
Crucial to the closed model structure is the operation $A\mapsto Cat(A)$ which
takes an $n$-precat to an $n$-category \cite{nCAT}. One may think of an
$n$-precat as a ``system of generators and relations'' for defining an
$n$-category, and of $Cat(A)$ as the $n$-category thusly defined.
This operation is explained in \cite{nCAT}, see also \cite{effective} and
\cite{descente}. To be brief, one can say that
$Cat(A)$ is obtained by throwing onto $A$ in a minimal way all of the
elements which are needed in order to satisfy the  definition of being an
$n$-category.  It is uniquely characterized (up to equivalence)
by the condition that there be a natural transformation $i_A:A\rightarrow
Cat(A)$
which (1) is an equivalence of $n$-categories whenever $A$ is already an
$n$-category, and (2) yields an equivalence of $n$-categories $Cat(i_A)$. (This
characterization is Proposition 4.2 of \cite{nCAT}, see also \cite{descente}).
\end{parag}

\begin{parag}
\label{othercat}
The other piece of information which we will need to know about the operation
$Cat$ is that it is partially composed of an operation $Fix$ which has the
effect
of operating $Cat$ on each of the $n-1$-categories $A_{p/}$; and that if
the Segal
maps $$
A_{p/}\rightarrow A_{1/}\times _{A_0} \ldots \times _{A_0} A_{1/}
$$
are weak equivalences of $n-1$-precats then $Cat(A)$ is equivalent to $Fix(A)$.

The other type of operation going into $Cat$ is called $Gen[m]$; we refer to
\cite{nCAT} for its description.
\end{parag}

\begin{parag}
\label{we}
{\em Weak equivalences:}
A morphism
$$
A\rightarrow B
$$
of $n$-precats is  a {\em weak equivalence} if the induced morphism of
$n$-categories
$$
Cat(A)\rightarrow Cat(B)
$$
is an equivalence of $n$-categories---cf \cite{Tamsamani} and
\ref{defequiv1}, \ref{essentialsurjectivity} above.
\end{parag}

\begin{parag}
\label{trivcofibs}
{\em Trivial cofibrations:}
A morphism $A\rightarrow B$ is said to be a {\em trivial cofibration}
if it is a cofibration and a weak equivalence.
\end{parag}

\begin{parag}
\label{fibs}
{\em Fibrations:}
A morphism $A\rightarrow B$ of $n$-precats is said to be a {\em fibration}
if it satisfies the following lifting property: for every trivial cofibration
$E'\hookrightarrow E$ and every morphism $E\rightarrow B$ provided with a
lifting over $E'$ to a morphism $E'\rightarrow A$, there exists an extension of
this to a lifting $E\rightarrow A$.

An $n$-precat $A$ is said to be {\em fibrant} if the canonical (unique)
morphism $A\rightarrow \ast$ to the constant presheaf with values one point,
is a fibration. 

A fibrant $n$-precat is, in particular, an $n$-category. This is because the
elements which need to exist to give an $n$-category may be obtained as
liftings of certain standard trivial cofibrations (those denoted $\Sigma
\rightarrow h$ in \cite{nCAT}).
\end{parag}

The following theorem allows us to do ``homotopy theory'' with $n$-categories.

\begin{theorem}
\label{cmc}
{\rm (\cite{nCAT} Theorem 3.1)}
The category $nPC$ of $n$-precats with the above classes of cofibrations, weak
equivalences and fibrations, is a closed model category.
\end{theorem}

\begin{parag}
\label{pushoutA}
We don't recall here the meaning of the closed model condition (axioms CM1-CM5
\cite{QuillenAnnals}). We just point out one of the main axioms in Quillen's
original point of view (Axiom M3 on page 1.1 of \cite{Quillen}) which says
that if
$A\rightarrow B$ is a trivial cofibration and $A\rightarrow C$ is any morphism
then $C\rightarrow B\cup ^AC$ is again a trivial cofibration.
This is a consequence of the axioms
CM1--CM5, but also one can remark that in the proof of \cite{nCAT}
(modelled on that of \cite{Jardine}) the main step which is done first
(\cite{nCAT} Lemma 3.2) is to
prove this property of preservation of trivial cofibrations by coproducts.
\end{parag}

\begin{parag}
\label{342}
Another of the axioms of a closed model category which we shall use heavily
is the ``three for the price of two'' axiom: given morphisms
$$
A\stackrel{f}{\rightarrow} B \stackrel{g}{\rightarrow}C
$$
in $nPC$, if any two of $f$, $g$ or $gf$ are weak equivalences then so is the
third.
\end{parag}

\bigskip

\subnumero{The construction $\Upsilon$}

The construction $\Upsilon$ introduced in \cite{limits} is one of the main tools
in the present note. The idea is that we would like to talk about the
basic $n$-category with two objects (denoted $0$ and $1$) and with a given
$n-1$-category $E$ of morphisms from $0$ to $1$ (but no morphisms in the other
direction and only identity endomorphisms of $0$ and $1$).  We call this
$\Upsilon (E)$.  To be more precise we do this on the level of precats: if
$E$ is
an $n-1$-precat then we obtain an $n$-precat $\Upsilon (E)$.  The main property
of this construction is that if $A$ is any $n$-category then a morphism of
$n$-precats  $$
f:\Upsilon (E) \rightarrow A
$$
corresponds exactly to a choice of two objects $x=f(0)$ and $y=f(1)$ together
with a morphism of $n-1$-precats $E\rightarrow A_{1/}(x,y)$.

One can see $\Upsilon (E)$ as the universal $n$-precat $A$ with two objects
$x,y$ and a map $E\rightarrow A_{1/}(x,y)$.

\begin{parag}
\label{upsilon}
We also need more general things of the form $\Upsilon ^2(E,F)$ having objects
$0,1,2$ and similarly a $\Upsilon ^3$.
(These will not have quite so simple an interpretation as
universal objects.)

Suppose $E_1,\ldots , E_k$ are $n-1$-precats. Then we define the $n$-precat
$$
\Upsilon ^k(E_1,\ldots , E_k)
$$
in the following way. Its object set is the set with $k+1$ elements denoted
$$
\Upsilon ^k(E_1,\ldots , E_k)_0 = \{ 0,\ldots , k\} .
$$
Set
$$
E_{ij}:= E_{i+1}\times \ldots \times E_{j-1} \times E_j
$$
with the convention that $E_{ii}=\ast$ and $E_{ij}=\emptyset$ if $j<i$.
Then define
$$
\Upsilon ^k(E_1,\ldots , E_k)_{p/}(y_0, \ldots , y_p)
:= E_{y_0y_1} \times \ldots \times E_{y_{p-1}y_p}.
$$
Notice that this is empty if any $y_i > y_j$ for $i<j$, equal to
$\ast$ if $y_0=\ldots = y_p$, and otherwise
$$
\Upsilon ^k(E_1,\ldots , E_k)_{p/}(y_0, \ldots , y_p):= E_{y_0 +1} \times \ldots
\times E_{y_p-1}\times E_{y_p}.
$$
\end{parag}

\begin{parag}
\label{erratum}
{\em Erratum:}
Note that the previous paragraph corrects an erratum in \cite{limits}
where the product in the definition of
$\Upsilon ^k(E_1,\ldots , E_k)_{p/}(y_0, \ldots ,
y_p)$ was erroneously written as starting with $E_{y_0}$ rather than
$E_{y_0 +1}$. Furthermore the notation from \cite{limits} erroneously suggested
that the product involved only  the $E_{y_i}$ rather than all the $E_i$
with $y_0<i\leq y_p$ as is the case.
\end{parag}

\begin{parag}
For example when $k=1$ (and we drop the superscript $k$ in this case)
$\Upsilon E$ is the $n$-precat with two objects $0,1$ and with $n-1$-precat of
morphisms from $0$ to $1$ equal to $E$. Similarly $\Upsilon ^2(E,F)$ has
objects $0,1,2$ and morphisms $E$ from $0$ to $1$, $F$ from $1$ to $2$ and
$E\times F$ from $0$ to $2$.  We picture $\Upsilon ^k(E_1,\ldots , E_k)$
as a $k$-simplex (an edge for $k=1$, a triangle for $k=2$, a tetrahedron for
$k=3$). The edges are labeled with single $E_i$, or products $E_i \times \ldots
\times E_j$.
\end{parag}

\begin{parag}
There are inclusions of these $\Upsilon^k$ according to the faces of the
$k$-simplex. The principal faces give inclusions
$$
\Upsilon ^{k-1}(E_1,\ldots , E_{k-1})\hookrightarrow \Upsilon ^k(E_1, \ldots ,
E_k),
$$
$$
\Upsilon ^{k-1}(E_2,\ldots , E_{k})\hookrightarrow \Upsilon ^k(E_1, \ldots ,
E_k),
$$
and
$$
\Upsilon ^{k-1}(E_1,\ldots , E_i\times E_{i+1}, \ldots  , E_{k})\hookrightarrow
\Upsilon ^k(E_1, \ldots , E_k).
$$
The inclusions of lower levels are deduced from these by induction. Note that
these faces $\Upsilon ^{k-1}$ intersect along appropriate $\Upsilon ^{k-2}$.
\end{parag}

\begin{remark}
\label{upsistar}
$\Upsilon (\ast )= I$ is the category with objects $0,1$ and with a unique
morphism from $0$ to $1$.  A map $\Upsilon (\ast )\rightarrow A$ is the same
thing as a pair of objects $x,y$ and a $1$-morphism from $x$ to $y$, i.e. an
object of $A_{1/}(x,y)$.
\end{remark}

\begin{parag}
\label{interpupsilon1}
See \cite{limits} 2.4.5-2.4.7 for
another way of constructing the $\Upsilon ^k$.
\end{parag}

\begin{parag}
\label{trivinclusions}
One thing which we often will need to know below is when an inclusion from a
union of faces, into the whole $\Upsilon ^k$, is a trivial cofibration.
For $k=2$ the inclusion which is a trivial cofibration is
$$
\Upsilon (E)\cup ^{\{ 1\} } \Upsilon (F) \hookrightarrow \Upsilon
^2(E, F).
$$

The reader is referred to \cite{limits} 2.4.8-2.4.10 for a similar
discussion for
$k=3$ (which  will not be used in the present note).
\end{parag}

\pagebreak[3]

\numero{Stabilization}

\subnumero{Minimal dimension}

\begin{parag}
In what follows we will be working in the closed model category $nPC$ of
$n$-precats (unless otherwise specified).
Recall that we have elements $1^i:= (1,\ldots , 1)\in \Theta ^n$ giving rise to
the representable $n$-precats $h(1^i)$ (i.e. the presheaves of sets
represented by $1^i$) \cite{nCAT}. Set
$$
F^i:=h(1^i)
$$
(for ``$i$-fl\^eche''). Let $\partial F^i$ denote the ``boundary''
i.e. the union of all of the lower degree arrows in $F^i$ (see \ref{interms}). A
map  $F^i\rightarrow A$ is the specification of an $i$-morphism in an $n$-precat
$A$. \end{parag}

\begin{parag}
One must include the limit case
$\partial F^{n+1} \rightarrow F^{n+1}$ where
$F^{n+1}:= F^n$ and
$$
\partial F^{n+1} := F^n\cup ^{\partial F^n} F^n.
$$
\end{parag}

\begin{parag}
\label{interms}
In terms of the notation $\Upsilon$, we have
$$
F^i= \Upsilon (F^{i-1}),
$$
starting with $F^0=\ast$, and also that
$$
\partial F^i = \Upsilon (\partial F^{i-1})
$$
starting with $\partial F^0=\emptyset$.
\end{parag}

\begin{parag}
\label{defmindim}
We say that a morphism of $n$-precats $f:A\rightarrow B$ has {\em minimal
dimension $m$} (usually denoted $m(f)$) if $m$ is the largest integer such that
there exists a diagram $$
\begin{array}{ccc}
A& \rightarrow & A' \\
\downarrow && \downarrow \\
B & \rightarrow & B'
\end{array}
$$
with the horizontal arrows being weak equivalences, the left vertical
arrow being the given one, and the right vertical arrow being a successive
(eventually transfinite) pushout by either trivial cofibrations or cofibrations
of the form  $\partial F^i \hookrightarrow F^i$ for $i\geq m$.

This means that homotopically, $B$ is obtained from $A$ by adding on
cells of dimension $\geq m$, and that $m$ is the largest such integer.
\end{parag}

\begin{parag}
\label{exists}
The minimal dimension as defined above exists, that is to say that there
always exists a diagram with the  right vertical arrow being a successive
pushout by either trivial cofibrations or cofibrations
of the form  $\partial F^i \hookrightarrow F^i$ for $i\geq 0$. To see this, note
that any cofibration of $n$-precats is for formal reasons a successive pushout
of cofibrations of the form $\partial h(M)\rightarrow h(M)$ for $M\in
\Theta ^n$.
If $M=(m,M')$ then  $$ \partial h(m,M')\hookrightarrow h(m,M')
$$
is weakly equivalent to a coproduct of inclusions of the form
$$
\partial h(1,M')\hookrightarrow h(1,M'),
$$
and we have
$$
h(1,M')=\Upsilon (h(M')), \;\;\;\; \partial h(1,M')=\Upsilon (\partial h(M')).
$$
Noting that $\Upsilon$ commutes with pushouts and preserves trivial
cofibrations,
and that by induction on $n$ we have the desired statement for $h(M')$, we get
that  $\partial h(1,M')\hookrightarrow h(1,M')$ admits a minimal dimension as in
\ref{defmindim}. Therefore $\partial h(M)\rightarrow h(M)$ admits a minimal
dimension, so any cofibration of $n$-precats admits a minimal dimension. By the
closed model structure any morphism is equivalent to a cofibration, so any
morphism admits a minimal dimension.
\end{parag}

\begin{parag}
\label{mequalsinf}
In the definition of minimal dimension \ref{defmindim} we make the
convention that
if $f$ is an equivalence then $m(f):= \infty$. Thus the minimal dimension of a
morphism of $n$-precats takes values in the set $\{ 0,1,\ldots , n,n+1, \infty
\}$. In particular if $m(f)\geq n+2$ then $f$ is an equivalence.
\end{parag}

\begin{parag}
\label{casenequals0}
For example if $n=0$ an $n$-precat is just a set and the minimal dimension of a
morphism $f$ of sets is $m(f)=\infty$ if $f$ is an isomorphism, $m(f)= 1$ if $f$
is surjective but not an isomorphism, and $m(f)=0$ otherwise.
\end{parag}

\begin{parag}
\label{property}
A consequence of the definition and the fact that the closed model category
$nPC$ is left proper (i.e. proper for cofibrant pushouts \cite{nCAT}
Theorem 6.7)
is the following property of minimal dimension: if $f:A\rightarrow B$ is a
morphism of minimal dimension $m(f)$  and if $g: A\rightarrow C$ is any
morphism,
and if either $f$ or $g$ is a cofibration, then the minimal dimension of the
morphism $$
C\rightarrow B\cup ^AC
$$
is at least $m(f)$.
\end{parag}

\subnumero{Relation with $k$-connectedness}

\begin{parag}
\label{kconnected}
We say that an $n$-category $A$ is {\em $k$-connected} (for $0\leq k\leq n$)
if the truncation \ref{highertau} $\tau _{\leq k}(A)$
is a contractible $k$-category (i.e. the morphism $\tau _{\leq k} (A)\rightarrow
\ast$ is an equivalence of $k$-categories).
\end{parag}

\begin{parag}
\label{kconnected2}
Say that an $n$-precat $A$ is {\em $k$-connected} if the $n$-category $Cat(A)$
is $k$-connected. The notion of $k$-connectedness is preserved by equivalence of
$n$-precats (the truncations are preserved by equivalences of
$n$-categories, see
\cite{Tamsamani}).
\end{parag}

\begin{parag}
\label{kconnected3}
If $A$ is $k$-connected for some
$k\geq 0$ then $\tau _{\leq 0}A$ is the one-point set, i.e. there
is a unique equivalence class of objects. Thus the choice of an object $a\in
A_0$ will be well-defined up to equivalence.
\end{parag}

\begin{lemma}
\label{kconnequivmindim}
Suppose $A$ is an $n$-category, and choose an object $a\in A_0$. Then $A$ is
$k$-connected if and only if the minimal dimension of the morphism $\{ a\}
\rightarrow A$ is $\geq k+1$.
\end{lemma}
{\em Proof:}
Suppose $A$ is a $k$-connected $n$-precat, and let
$$
B:= A\cup ^{\partial F^i}F^i
$$
with $i\geq k+1$. We will show that $B$ is also $k$-connected. Note that
$$
Cat(B) \cong Cat(Cat(A)\cup ^{\partial F^i}F^i)
$$
so we may assume that $A$ is an $n$-category. Fix objects $x_0,\ldots ,
x_p\in A$.
Note that
$A_{p/}(x_0,\ldots , x_p)$ is $k-1$-connected, in fact
$$
\tau _{\leq k}(A)_{p/}(x_0,\ldots ,x_p) =
\tau _{\leq k-1}(A_{p/}(x_0,\ldots , x_p)
$$
(this follows from Tamsamani's definition of the truncation operations
``starting from the top and going down'').
Now the morphism of $n-1$-precats
$$
(\partial F^i)_{p/}(y_0,\ldots , y_p)\rightarrow
F^i_{p/}(y_0,\ldots , y_p)
$$
has minimal dimension $\geq i-1$ (this can be seen by the expressions
$F^i=\Upsilon (F^{i-1})$ and $\partial F^i=\Upsilon (\partial F^{i-1})$
together with the definition of $\Upsilon$ and Theorem \ref{mainth} for
$n-1$-categories).
Therefore by the inductive version of the statement of the present lemma for
$n-1$-precats, we get that
$$
(A\cup ^{\partial F^i}F_i)_{p/}(x_0,\ldots , x_p)
$$
is $k-1$-connected. Finally, recall that $Cat$ is, up to equivalence, a
successive composition of operations denoted $Fix$ and $Gen[m]$. A
pushout of $k-1$-connected $n-1$-precats remains $k-1$-connected, so the
operation
$Gen[m]$ (cf \cite{nCAT}) preserves the property  of $k-1$-connectedness of the
components $(\; )_{p/}(x_0,\ldots , x_p)$; the operation $Fix$ also clearly
does.
Therefore
$$
Cat(A\cup ^{\partial F^i}F_i)_{p/}(x_0,\ldots , x_p)
$$
is $k-1$-connected, which implies that $Cat(A\cup ^{\partial F^i}F_i)$
is $k$-connected.

For the purposes of the above argument when $k=0$
we use ``$-1$-connected'' to mean nonempty.

We have now obtained one half of the lemma:  if the minimal dimension
of the morphism $\{ a\} \rightarrow A$ is $\geq k+1$ then $A$ is $k$-connected.

\begin{parag}
\label{prooflemma2}
The proof of the other half of the lemma involves a construction which is of
some interest regarding the definitions of Baez-Dolan \cite{BaezDolan}---we
replace a $k$-connected $n$-category by an equivalent $n$-category which is
trivial (i.e. {\em equal} to $\ast$) in degrees $\leq k$. This is related to
Dunn's Segal-type $E_k$-machine, cf the discussion in \cite{descente}.

Suppose that $A$ is an $n$-category. Choose an object
$a\in A$. Let $A'\subset A$ be the $n$-precat defined by setting (for $M\in
\Theta ^n$)
$$
A'_{M}\subset A_{M}
$$
equal to the subset of elements $\alpha $ such that for any morphism
$u:U\rightarrow M$ in $\Theta ^n$ with $| U | \leq k$, the image
$u^{\ast}(\alpha )$ is equal to the degeneracy $d^{\ast}(a)$ where
$d: U\rightarrow 0$ is the unique map in $\Theta ^n$.

By construction we have $A'_M =\ast$ whenever $|M|\leq k$.
\end{parag}

\begin{parag}
\label{nameit}
We introduce a notation for the construction defined above:
$$
{\bf Wh}_{> k}(A,a):= A'.
$$
The letters $Wh$ refer to the fact that in the context of spaces, this
is the {\em Whitehead tower}.
\end{parag}

\begin{parag}
\label{prooflemma3}
We claim that if $A$ is an $n$-category, $a\in A_0$ is an object, and $0\leq k
\leq n$, then ${\bf Wh}_{> k}(A,a)$ is again an $n$-category. Concurrently we
claim that  if
$f:(A,a)\rightarrow (B,b)$ is an equivalence sending $a$ to $b$ then  the
induced
morphism ${\bf Wh}_{> k}(f)$ is an equivalence.

To
prove these claims (by induction on $n$), note first that for $k=0$ we are just
taking the full sub-$n$-category of $A$ containing only the object $a$, which is
again an $n$-category. The statement about equivalences is equally clear in this
case. If $k\geq 1$, note that
$$
{\bf Wh}_{> k}(A,a)_{p/} = {\bf Wh}_{> (k-1)}(A_{p/}(a,
\ldots , a),d_p(a))
$$
where $d_p(a)$ is the degeneracy of $a$ considered as an object of
$A_{p/}(a, \ldots , a)$. By the first claim for $n-1$ the components
${\bf Wh}_{> k}(A,a)_{p/}$ are $n-1$-categories.
The fact
that ${\bf Wh}_{> k}$ is compatible with direct products,
together with the second of our claims for $n-1$, imply that the Segal maps
for the
${\bf Wh}_{> k}(A,a)_{p/}$ are equivalences. This proves that
${\bf Wh}_{> k}(A,a)$ is an $n$-category. The second claim about equivalences
follows from the fact that ${\bf Wh}_{> k}(f)$ is evidently essentially
surjective, and full faithfulness comes (via the above formula for
${\bf Wh}_{> k}(A,a)_{1/}(a,a)$) from the inductive statement for $n-1$. This
completes the proof of the two claims.
\end{parag}

\begin{parag}
\label{prooflemma3bis}
We claim furthermore that if $A$ is $k$-connected then the morphism
$$
i:{\bf Wh}_{> k}(A,a)\rightarrow A
$$
is an equivalence of $n$-categories. For any $k\geq 0$,  $k$-connectedness
implies that $\tau _{\leq 0}(A)$ has only one object, so the morphism $i$ is
essentially surjective. In the case $k=0$, ${\bf Wh}_{>
0}(A,a)_{1/}(a,a)=A_{1/}(a,a)$ so $i$ is fully faithful; this treats the case
$k=0$. For  $k\geq 1$ we proceed by induction on $n$ and use the same formula as
in \ref{prooflemma3} to get that $i$ is fully faithful.
\end{parag}

\begin{parag}
\label{prooflemma4}
Return to the proof of Lemma \ref{kconnequivmindim} where we suppose that
$A$ is $k$-connected. Set $A'={\bf Wh}_{> k}(A,a)$.
The morphism
$\{ a\} \rightarrow A'$ is obtained by a sequence of cofibrations of the form
$\partial h(M)\rightarrow h(M)$ for $| M| \geq k+1$.
By the same argument as in \ref{exists}, this morphism is obtained up to weak
equivalence by adding on cofibrations of the form $\partial F^i\hookrightarrow
F^i$ for $i\geq k+1$. Thus the minimal dimension of $\{ a\} \rightarrow A'$ is
at least $k+1$, hence (by \ref{prooflemma3bis}) the same is true of $\{ a\}
\rightarrow A$. This completes the proof of the lemma.
\eop
\end{parag}

\begin{parag}
\label{relnwithtruncation}
In the above discussion we have used in an essential way that the source of the
map in question consists of only one point.  For a general map, being of minimal
dimension $\geq k$ is related to the behavior on truncations but not exactly
the same thing. Suppose $f: A\rightarrow B$ is a morphism of $n$-categories.
If $m(f) \geq k+1$ then $\tau _{\leq k-1}(f)$ is an equivalence; and in the
other
direction, if $\tau _{\leq k}(f)$ is an equivalence then $m(f)\geq k+1$. One
cannot in general make a more precise statement than that.
\end{parag}

\begin{parag}
\label{relnwithtruncation2}
The behavior in the previous paragraph may be understood by reference to
classical homotopy theory (which is the case where the $n$-categories in
question are $n$-groupoids cf \cite{Tamsamani}). In this case, $f:A\rightarrow
B$ has minimal dimension $\geq k+1$ if and only if $f$ induces an isomorphism on
$\pi _i$ for $i\leq k-1$, and a surjection on $\pi _{k}$---minimal dimension
here corresponds to the dimensions of the cells which $f$ adds to $A$ to obtain
$B$.  Of course, if $A$ consists of only one point, then $f$ being surjective on
$\pi _{k}$ is the same thing as $\pi _{k}(f)$ being an isomorphism
and we recover the characterization of \ref{kconnequivmindim}.
\end{parag}

\subnumero{The main estimate}
Our main theorem is the following:
\begin{theorem}
\label{mainth}
Suppose $f:A\rightarrow B$ and $g:C\rightarrow D$ are cofibrations of
$n$-precats with minimal dimensions $m(f)$ and $m(g)$ respectively. Then the
minimal dimension of the morphism $$
f\wedge g: A\times D \cup ^{A\times C} B\times C \rightarrow B\times D
$$
is greater than or equal to $m(f) + m(g)$.
\end{theorem}

\begin{parag}
This is the analogue for $n$-categories of the visually obvious corresponding
statement in topology.
\end{parag}

\begin{parag}
\label{casezero}
We first indicate the proof for $n=0$ i.e. for sets.
It suffices to consider the morphisms
$$
a: \partial F^0 = \emptyset \rightarrow \ast = F^0,
$$
and
$$
b: \partial F^1 = 2\ast \rightarrow \ast = F^1,
$$
where $2\ast$ is the set with two elements. Note that $m(a)=0$ and $m(b)=1$.
Obviously $m(a\wedge a)=0$. We have
$$
a\wedge b = \left( \emptyset \cup ^{\emptyset} 2\ast \rightarrow \ast \right)
$$
so $m(a\wedge b)=1$, and similarly $m(b\wedge a)=1$. Finally
$$
b\wedge b = \left( {2\ast} \cup ^{2\ast \times 2\ast } {2\ast} \rightarrow \ast
\right)
$$
where the two maps in the coproduct are the two projections; thus
$b\wedge b$ is an isomorphism so $m(b\wedge b)=\infty \geq 1+1$.
This completes the verification of Theorem \ref{mainth} for $n=0$.
\end{parag}

\begin{parag}
\label{dirprod}
Note that the direct product with a weak equivalence is again a weak equivalence
(see \cite{nCAT} Theorem 5.1). This treats the case of \ref{mainth} where
one of $f$ or $g$ is an equivalence (i.e. has minimal dimension $\infty$).
\end{parag}

\begin{parag}
We now give the proof of the theorem by induction on $n\geq 1$.
We suppose in what follows that the theorem is known for
$n-1$-categories. The main lemma is the following.
\end{parag}

\begin{lemma}
\label{thelemma}
Suppose $f:A\rightarrow B$ and $g:C\rightarrow D$ are cofibrations of
$n-1$-precats with minimal dimensions $m(f)$ and $m(g)$ respectively. Then the
minimal dimension of
$$
\Upsilon (A) \times \Upsilon (D)
\cup ^{\Upsilon (A) \times \Upsilon (C)}
\Upsilon (B) \times \Upsilon (C)
$$
$$
\rightarrow
\Upsilon (B) \times \Upsilon (D)
$$
is at least $m(f) + m(g) + 2$.
\end{lemma}

\begin{parag}
Using \ref{dirprod} and the fact that $F^i= \Upsilon (F^{i-1})$ and
$\partial F^i = \Upsilon (\partial F^{i-1})$ it is easy to deduce the theorem
from the lemma. So, we will prove the lemma.
\end{parag}

\begin{parag}
\label{mdupsilon}
Note also that using \ref{interms} and \ref{defmindim} one obtains that for any
map $h: X\rightarrow X'$ of $n-1$-precats of minimal dimension $m(h)$, the
minimal dimension of  $$
\Upsilon (h): \Upsilon (X)\rightarrow \Upsilon (X')
$$
is at least equal to $m(h)+1$.
\end{parag}

\begin{parag}
\label{square}
We can write
$$
\Upsilon (B)\times \Upsilon (D) =
$$
$$
\Upsilon ^2(B, D) \cup ^{\Upsilon (B\times D)} \Upsilon ^2(D,B)
$$
(the reader is urged to draw a square divided into two triangles, labeling the
horizontal edges with $B$, the vertical edges with $D$, and the hypotenuse with
$B\times D$). The maps in the coproduct are cofibrations. The
equation can be checked directly from the definition of $\Upsilon$.

We have similar equations for the other products occuring in the statement
of the lemma.
\end{parag}

\begin{parag}
\label{rem1}
We now note that
\end{parag}
$$
\left( \Upsilon (A) \cup ^{\ast} \Upsilon (D) \right)
\cup ^{(\Upsilon (A) \cup ^{\ast} \Upsilon (C))}
\left( \Upsilon (B) \cup ^{\ast} \Upsilon (C) \right)
$$
$$
= \Upsilon (B)\cup ^{\ast} \Upsilon (D).
$$
In this formula the coproducts are taken over the $1$-object sets $\ast$, for
example in $\Upsilon (A)\cup ^{\ast} \Upsilon (D)$, $1\in \Upsilon (A)$ is
joined with $0\in \Upsilon (D)$.

\begin{parag}
\label{proper}
The closed model category $nPC$ is left proper (see \cite{nCAT} Theorem 6.7),
which implies that a triple of equivalences on each of the elements of a
cofibrant
pushout diagram yields an equivalence of the pushout. This consequence of
properness is known as ``Reedy's lemma'' \cite{Reedy} and can be found in any
number of recent references on closed model category theory.
\end{parag}

\begin{parag}
\label{notnQ}
Introduce the notation
$$
Q(A,B,C,D):= \Upsilon ^2(A, D) \cup ^{\Upsilon ^2(A,C)}
\Upsilon ^2(B,C).
$$
The morphisms
$$
i(A,D): \Upsilon (A)\cup ^{\ast} \Upsilon (D) \rightarrow
\Upsilon ^2(A,D),
$$
and similarly $i(A,C)$ and $i(B,C)$
are all equivalences (\ref{trivinclusions}). Using
the equation \ref{rem1} we obtain a morphism
$$
\Upsilon (B)\cup ^{\ast} \Upsilon (D) \rightarrow Q(A,B,C,D),
$$
and by the remark \ref{proper} this morphism is an equivalence.
Then using the fact that $i(B,D)$ is an equivalence and the
three-for-the-price-of-two property \ref{342} we get that
the morphism
$$
u: Q(A,B,C,D)\rightarrow
\Upsilon ^2(B,D)
$$
is an equivalence.
Similarly the map
$$
v: Q(C,D,A,B)= \Upsilon ^2(D, A) \cup ^{\Upsilon ^2(C,A)}
\Upsilon ^2(C,B) \rightarrow \Upsilon ^2(D,B)
$$
is an equivalence.
\end{parag}

\begin{parag}
\label{rem2}
Using the equations \ref{square} and a manipulation of pushout formulae,
we obtain
$$
\Upsilon (A)\times \Upsilon (D) \cup ^{\Upsilon (A)\times \Upsilon (C)}
\Upsilon ( B) \times \Upsilon (C)
$$
$$
=
$$
$$
Q(A,B,C,D) \cup ^Y  Q(C,D, A, B)
$$
where
$$
Y:= \Upsilon (A\times D \cup ^{A\times C} B\times C).
$$
In particular the left side of the morphism in the lemma is
$Q(A,B,C,D) \cup ^Y  Q(C,D, A, B)$.
\end{parag}

\begin{parag}
\label{rem3}
Now the equivalences $u$ and $v$ of \ref{notnQ} give an equivalence
$$
Q(A,B,C,D) \cup ^Y  Q(C,D, A, B)\rightarrow
\Upsilon ^2(B,D) \cup ^{Y} \Upsilon ^2(D,B).
$$
Rewrite the right side of this morphism as
$$
V:= \Upsilon ^2(B,D) \cup ^{Y} \Upsilon ^2(D,B)=
\Upsilon ^2(B,D) \cup ^{\Upsilon (B\times D)} Z
\cup ^{\Upsilon (D\times B)} \Upsilon ^2(D,B)
$$
where
$$
Z:=\Upsilon (B\times D) \cup ^{Y} \Upsilon (B\times D).
$$
(The second map $\Upsilon (D\times B)\rightarrow Z$ in the expression for $V$
is obtained using the standard isomorphism $D\times B \cong B\times D$. )

For the lemma, it suffices to show that the morphism
$$
w: \Upsilon ^2(B,D) \cup ^{Y} \Upsilon ^2(D,B)
\rightarrow
\Upsilon (B)\times \Upsilon (D),
$$
has minimal dimension as required.
\end{parag}

\begin{parag}
\label{rem4}
We have a morphism
$$
\alpha : Z\rightarrow \Upsilon (B\times D)
$$
(the identity on both components).
Taking the coproduct of the expression $V$ from \ref{rem3} along the morphism
$\alpha$ gives
$$
V\cup ^Z \Upsilon (B\times D) = \Upsilon ^2(B,D) \cup ^{\Upsilon (B\times D)}
\Upsilon ^2(D,B)
$$
$$
= \Upsilon (B)\times \Upsilon (D).
$$
In other words, the morphism $w$ is obtained by a coproduct of a cofibration
$Z\hookrightarrow V$ along $\alpha$. In particular
the minimal dimension of the arrow in the lemma is at least as big as the
minimal dimension $m(\alpha )$ which we shall now bound.
\end{parag}

\begin{parag}
\label{rem5}
By the global inductive hypothesis of our proof of the theorem, we may assume
that Theorem \ref{mainth} is known for $A$, $B$, $C$, $D$ which are
$n-1$-categories. Therefore, setting
$$
W:= A\times D \cup ^{A\times C} B\times C,
$$
the minimal dimension of the cofibration
$$
W\rightarrow B\times D
$$
is at least $m(f)+m(g)$.
\end{parag}

\begin{parag}
\label{rem6}
Note with this notation that $Y=\Upsilon (W)$, and
$Z= \Upsilon (B\times D \cup ^W B\times D)$. The morphism $\alpha$ is obtained
by
applying $\Upsilon$ to the map
$$
\beta : B\times D \cup ^W B\times D \rightarrow B\times D.
$$
\end{parag}

\begin{parag}
\label{claim1}
\, {\bf Claim:} if $i: E\hookrightarrow F$ is a cofibration of $n-1$-precats of
minimal dimension $m(i)$ then the minimal dimension of the morphism
$$
F\cup ^EF\rightarrow F
$$
is at least $m(i)+1$.
\end{parag}

\begin{parag}
\label{completeproof}
Using this claim, we can complete the proof of Lemma \ref{thelemma}.
Apply the claim to the cofibration $W\rightarrow B\times D$ of \ref{rem5}
which has minimal dimension at least $m(f)+m(g)$. By the claim we find that the
minimal dimension of the map $\beta$ in \ref{rem6} is at least $m(f)+m(g)+1$.
Since $\alpha =\Upsilon (\beta )$, applying \ref{mdupsilon} we get that
$m(\alpha )\geq m(f)+m(g)+2$. In view of \ref{rem4} this will prove the lemma.
\end{parag}

\begin{parag}
\label{proveclaim}
To finish, we prove Claim \ref{claim1}. Note that it concerns $n-1$-precats
so by the
global inductive hypothesis we may use the statement of Theorem \ref{mainth}.
Also, note
that the
statement of the claim is easily verified by hand for $n-1=0$ i.e. when $E$ and
$F$ are sets. Therefore we may assume here that $n\geq 2$, i.e. $E$ and $F$ are
$n-1$-precats with $n-1\geq 1$.

Let $\overline{I}$ be the $1$-category with two
objects $0,1$ and a single isomorphism between them; we shall consider it as an
$n-1$-precat.
The minimal dimension of the cofibration of $n-1$-precats
$$
j:\{ 0,1\} \rightarrow \overline{I}
$$
is $1$. Therefore, applying the inductive hypothesis that Theorem \ref{mainth}
holds for $n-1$-precats, we get that the minimal dimension of the map
$$
i\wedge j: E\times \overline{I} \cup ^{E\times \{ 0,1\} }
F\times \{ 0,1\}
\rightarrow F \times \overline{I}
$$
is at least $m(i)+1$. Finally, note that
$$
E\times \overline{I} \cup ^{E\times \{ 0,1\} }
F\times \{ 0,1\}  =
$$
$$
\left( F\times \{ 0\} \cup ^{E\times \{ 0\}}E\times \overline{I} \right)
\cup ^{E\times \overline{I}}
\left( F\times \{ 1\} \cup ^{E\times \{ 1\}}E\times \overline{I} \right) ,
$$
and this latter expression maps to the coproduct
$F\cup ^EF$ by a map which comes from the three equivalences
$$
\left( F\times \{ 0\} \cup ^{E\times \{ 0\}}E\times \overline{I} \right)
\rightarrow F,
$$
$$
E\times \overline{I} \rightarrow E,
$$
and
$$
\left( F\times \{ 1\} \cup ^{E\times \{ 1\}}E\times \overline{I} \right)
\rightarrow F.
$$
(Use the results of \cite{nCAT}---preservation of weak equivalences by direct
product, and preservation of trivial cofibrations by pushout---to prove that
these maps are equivalences).
As recalled at \ref{proper}, the closed model category $nPC$ is left-proper so a
map of pushout diagrams composed of three weak equivalences, induces a weak
equivalence on pushouts. Therefore we have a diagram
$$
\begin{array}{ccc}
E\times \overline{I} \cup ^{E\times \{ 0,1\} }
F\times \{ 0,1\} & \rightarrow & F\times \overline{I} \\
\downarrow && \downarrow \\
F\cup ^EF & \rightarrow & F
\end{array}
$$
in which the vertical arrows are weak equivalences, and in which the top map has
(as we have seen previously) minimal dimension at least $m(i)+1$. Therefore the
bottom map has minimal dimension at least $m(i)+1$, which is Claim \ref{claim1}.
This completes the proof of Lemma \ref{thelemma}.
\eop
\end{parag}

\begin{parag}
\label{completetheorem}
The theorem follows from the lemma using the definition and the fact that the
basic cells $F^i$ may be considered as arising from $\ast$ by iterated
applications of the operation $\Upsilon$ cf \ref{interms}.
\eop
\end{parag}

\subnumero{The stabilization hypothesis}
We now show how to prove Theorem \ref{thmintro} of the introduction.

\begin{parag}
\label{basepointed}
By a {\em pointed $n$-precat} we mean an $n$-precat $A$ with chosen object $a$
or equivalently with a morphism
$a: \ast
\rightarrow A$. The {\em minimal dimension} of a pointed $n$-precat is defined
as the minimal dimension of the map $a$.
\end{parag}

\begin{parag}
\label{wh}
Recall from \ref{prooflemma2}, \ref{nameit}, \ref{prooflemma3} that we have
defined the {\em Whitehead operation}  $(A,a) \mapsto {\bf Wh}_{>k}(A,a)$ which
gives a subobject of $A$ containing as $i$-morphisms only the $1^i_a$ for $i\leq
k$. Recall also \ref{prooflemma3bis} that if $A$ is $k$-connected then the
morphism   $$
A\rightarrow {\bf Wh}_{>k}(A,a)
$$
is an equivalence.
\end{parag}

\begin{parag}
\label{relnBD}
If $(A,a)$ is a pointed $n+k$-category
of minimal dimension $k$
then $A$ is $k-1$-connected \ref{kconnequivmindim} and equivalent to its
$n+k$-subcategory ${\bf Wh}_{>k-1}(A,a)$. This latter, in Baez-Dolan's
terminology, is a $k$-uply monoidal $n$-category. Conversely it is clear that  a
$k$-uply monoidal $n$-category, considered as a pointed $n+k$-category, has
minimal dimension $k$.

In what follows we trade $n+k$ for $n$ and shall look at pointed
$n$-categories $(A,a)$ which are $k-1$-connected or equivalently, of minimal
dimension $k$.
\end{parag}

\begin{corollary}
\label{estimate2}
Suppose $(A,a)$ and $(B,b)$ are pointed $n$-precats with minimal dimensions
$m(A,a)$ and $m(B,b)$ respectively. Then the minimal dimension of the morphism
$$
A\cup ^{\ast} B \rightarrow A\times B
$$
is at least $m(A,a) + m(B,b)$. In particular if $m(A,a)+m(B,b)\geq n+2$ then
the above morphism is a weak equivalence. (The coproduct in the display
identifies
$\ast = \{ a\}\subset A$ with $\ast = \{ b\} \subset B$.)
\end{corollary}
{\em Proof:}
Apply Theorem \ref{mainth} to the inclusions $\{ a\} \hookrightarrow A$ and
$\{ b\} \hookrightarrow B$. The map of Theorem \ref{mainth} for these two
cofibrations is exactly the map $A\cup ^{\ast} B \rightarrow A\times B$.
For the second statement, recall from \ref{mequalsinf} that a morphism of
minimal dimension $\geq n+2$ in fact has minimal dimension $\infty$ and is
a weak
equivalence.
\eop

\begin{parag}
\label{construction}
Now suppose $(A,a)$ is a pointed $n$-precat. We define a simplicial $n$-precat
$X$ by the following:
$$
X_{p/}:= A\cup ^{\ast} A \cup ^{\ast} \ldots \cup ^{\ast} A
\;\;\;\; (k \; \mbox{times}).
$$
In particular $X_{0/}$ (which we denote $X_0$) is the set $\ast$. The maps
in the
simplicial structure are obtained by using the identity map or else the
projection $p:A\rightarrow \ast$ on each component in an appropriately organized
way. For example
$$
X_{2/} = A\cup ^{\ast} A,
$$
the first face map is the identity on the first component and the projection
$p$ on the second component; the second face map is the identity on the second
component and the projection on the first component; and finally the $02$ face
map is the identity on both components.
\end{parag}

\begin{parag}
\label{Xprecat}
According to \ref{precat2},
we can (and will) consider the simplicial object  $X: \Delta ^o \rightarrow nPC$
as being an $n+1$-precat. Let  $Cat(X)$ denote the replacement of $X$
(considered
as an $n+1$-precat) by an
$n+1$-category \ref{opcat}.
\end{parag}

\begin{corollary}
\label{prethm1}
Suppose $(A,a)$ is a  pointed  $n$-category which has minimal dimension $k$.
Suppose that $2k \geq n+2$. Then the simplicial $n$-precat $X$ defined above
satisfies the Segal condition that the maps
$$
X_{p/} \rightarrow X_{1/} \times _{X_0} \ldots \times _{X_0} X_{1/}
$$
are weak equivalences. The morphisms
$$
X_{p/} \rightarrow Cat(X)_{p/}
$$
are weak equivalences. In particular we have the equivalence of $n$-categories
$$
A\stackrel{\cong}{\rightarrow} Cat(X)_{1/} (x , x ),
$$
where $x$ is the unique object of $X$, and $(Cat(X),x)$ has minimal dimension
$k+1$.
\end{corollary}
{\em Proof:} For the first part it follows from the previous corollary, noting
that $X_0=\ast$. The second part follows from the fact that the operation
$Cat$ in \cite{nCAT} may be viewed as starting with an operation $Fix$ which has
the effect of doing $Cat$ on each of the components $X_{p/}$. Under the
condition of the first part of the corollary that the Segal maps are weak
equivalences of $n$-precats, we get that the Segal maps for $Fix(X)$ are
equivalences of $n$-categories, i.e. $Fix(X)$ is already an $n+1$-ctaegory. In
this case $Fix(X)\rightarrow Cat(X)$ is an equivalence (the remaining operations
$Gen[m]$ that go into $Cat$ don't have any effect if the input is already an
$n+1$-category).  \eop

The following corollary is Theorem \ref{thmintro} of the introduction.

\begin{corollary}
\label{thm1}
Suppose $A$ is a $k$-connected $n$-category with $2k\geq n$. Then there exists
a $k+1$-connected $n+1$-category $Y$ with $Y_{1/}(y,y) \cong A$ (for the
essentially unique object $y$ of $Y$).
\end{corollary}
{\em Proof:} Set $Y:= Cat(X)$ in the previous corollary. Use Lemma
\ref{kconnequivmindim} to compare $k$-connectedness with minimal dimension.
\eop

\begin{parag}
\label{delooping}
This corollary says that $Y$ is a ``delooping'' of $A$. That basically gives the
Stabilization Hypothesis of \cite{BaezDolan} \cite{cat}, although there are
a lot
of further details which need to be considered: one would like to give an
equivalence between the basepointed $k$-connected $n$-categories and the
basepointed $k+1$-connected $n+1$-categories.
\end{parag}

\begin{parag}
\label{numerology}
{\em Numerology:} In
Baez-Dolan's notation a ``$k$-uply monoidal $n$-category'' is
a $k-1$-connected $n+k$-category cf \ref{relnBD}. Thus, the corollary
says that if $A$ is a $k$-uply monoidal $n$-category, and if  $2(k-1) \geq n+k$
(i.e. $k\geq n+2$), then there exists a delooping $Y$ of $A$.

\numero{Constructing examples}

In this section we will discuss some ways of constructing examples of
$k$-uply monoidal $n-k$-categories (i.e. $k-1$-connected $n$-categories).
Almost all of these examples appear in Baez-Dolan \cite{BaezDolan}
\cite{cat} (and in turn they are often due, for  low values of $k$ and $n$, to
other people before that---see the reference lists of \cite{BaezDolan},
\cite{cat}, \cite{BreenRecent}). In these cases our only contribution is to
clarify how they fit into the picture of Tamsamani's definition of
$n$-categories. In the last subsection on ``cohomological twisting'' we
present a
construction which might be new.

\subnumero{Monoidal and symmetric monoidal objects}

\newparag{monoidal1}
If $C$ is a category admitting direct products, recall that a {\em monoidal
object} in $C$ is an object $A$ together with a composition law
$m:A\times A\rightarrow A$, and a  morphism $e: \ast \rightarrow A$,
such that $m$ is associative and $e$ is an identity (these conditions are
expressed in terms of commutative diagrams). A monoidal object is {\em symmetric
monoidal} if $m$ satisfies the rule $m\circ \sigma = m$ where $\sigma :
A\times A\rightarrow A\times A$ is the morphism which interchanges the two
factors.

\newparag{monoidal2}
If $(A, m, e)$ is a monoidal object in $C$, then it determines in a natural way
a category $c(A,m,e)$ enriched over $C$, with only one object which we
denote $x$,
and
$$
Hom_{c(A,m,e)}(x,x)=A
$$
with $m$ as composition and $e$ as unit.

\newparag{monoidal3}
If $(A,m,e)$ is a symmetric monoidal object, then $c(A, m, e)$ itself is
naturally
endowed with a structure of symmetric monoidal object ($m$ again serves as the
multiplication), so one can iterate this construction to obtain a strict
$k$-category $c^k(A, m, e)$ with top-level morphisms enriched over $C$. To be
precise, this is a strict $k$-category with only one $i$-morphism $1^i_x$
for any
$i<k$ (we call the base object $x= 1^0_x$), and such that the $k$-endomorphisms
of $1^{k-1}_x$ are the object $A\in C$ with $m$ as multiplication. (This
structure may be taken as the definition of what we mean by ``strict
$k$-category
with top-level morphisms enriched over $C$'', see below for another
interpretation).

\newparag{monoidal4}
Conversely, any strict $k$-category with $k=1$ (resp. $k\geq 2$) with top-level
morphisms enriched over $C$, and having only one $i$-morphism $1^i_x$ for $i<k$,
is of the form $c^k(A,m,e)$ for a monoidal object (resp. symmetric monoidal
object) $(A,m,e)$ in $C$.

\newparag{multisimplicialinterp}
We can give a multisimplicial interpretation of what is meant by ``strict
$k$-category with top-level morphisms enriched over $C$''. Let $C^{\sqcup}$
denote the category obtained by formally adjoining disjoint sums to $C$
(in general we will use these
remarks in the case where $C$ admits disjoint sums, thus $C\cong C^{\sqcup}$
and the reader can ignore this notation if he wishes).  We obtain in
particular a
functor $Set \rightarrow C^{\sqcup}$ sending a set $S$ to the formal
disjoint sum
of $S$ copies of the final object $\ast\in C$ (this final object exists
because we
assume that $C$ admits products). We assume for example that this functor is
fully faithful. A {\em strict $k$-category with top-level morphisms
enriched over
$C^{\sqcup}$} is a functor
$$
A: (\Theta ^k)^o \rightarrow C^{\sqcup}
$$
such that the Segal maps (cf \ref{ncatsdef1}) are isomorphisms, and such
that the
images of $M\in \Theta ^k$ with $|M| <k$ are sets (i.e. contained in the
subcategory $Set\subset C^{\sqcup}$).

We adopt the same conventions as in \ref{arrows} for the set
of $i$-arrows in $A$, $i<k$. There is a $C^{\sqcup}$-object of $k$-arrows
$A_{1,\ldots , 1}$, and for any two $k-1$-arrows $u,v$ having $s(u)=s(v)$ and
$t(u)=t(v)$ we get $A_{1,\ldots , 1}(u,v)\in C^{\sqcup}$. If these latter
objects are actually in $C$ then we say that $A$ is a
{\em strict $k$-category with top-level
morphisms enriched over $C$}.

In the case where we require, furthermore, that the only $i$-arrows for $i<k$ be
the $1^i_x$, then we don't need to use $C^{\sqcup}$ in the above definition:
such an object $A$ is just a functor $\Theta ^k\rightarrow C$ such that
$M\mapsto
\ast \in C$ if $|M|<k$, and such that the Segal maps are isomorphisms.
(Compare Dunn \cite{Dunn}).

We leave to the reader the verification of the content of
\ref{monoidal2}, \ref{monoidal3} and \ref{monoidal4} with these precise
definitions.

\newparag{apply}
We now apply the above discussion to the category $C= nPC$ of $n$-precats.
Note that sums exist so $C\cong C^{\sqcup}$.
We obtain a notion of monoidal (resp. symmetric monoidal) object in $nPC$,
and such an object corresponds (with $k=1$ resp. $k\geq 2$) to a functor
$$
c^k(A):(\Theta ^k)^o\rightarrow nPC,
$$
such that the Segal maps (at all levels cf \cite{Tamsamani}) are isomorphisms,
and such that $c^k(A)$ takes on values which are sets on
$M\in \Theta ^k$ with $|M|
<k$. By definition then, $c^k(A)$ is an $n+k$-precat $A\in (n+k)PC$.
If the image lies in the subcategory $n-Cat \subset nPC$ of $n$-categories,
then $c^k(A)$ is an $n+k$-category (the Segal maps below level $k$
being isomorphisms
and above level $k$ being equivalences).

In general, let $Cat_{(n)}$ (resp. $Cat_{(n+k)}$) be the operations
on $nPC$ (resp. $(n+k)PC$)
yielding $n$-categories (resp. $n+k$-categories) \ref{opcat} cf \cite{nCAT}.
Then
we obtain a functor
$$
Cat _{(n)}\circ A : (\Theta ^k)^o\rightarrow n-Cat.
$$
The operation $Cat_{(n)}$ is not compatible with direct products but up to
equivalence it is (\cite{nCAT} Theorem 5.1). Therefore the Segal maps for
$Cat _{(n)}\circ c^k(A)$ remain equivalences and
$Cat _{(n)}\circ c^k(A)$ is an $n+k$-category.
Note that
$$
(Cat _{(n)}\circ c^k(A))_{1,\ldots , 1}(1^{k-1}_x, 1^{k-1}_x) =
Cat _{(n)}(c^k(A)_{1,\ldots  , 1}
(1^{k-1}_x, 1^{k-1}_x) =Cat _{(n)}(A)
$$
and this is weak equivalent to $A$. Thus we can say that $Cat _{(n)}\circ
c^k(A)$
is a $k$-uply monoidal $n$-category (or $k-1$-connected $n+k$-category)
with underlying $n$-category weakly equivalent to $A$.

By an inductive
application $k$ times of the statement that the operation $Cat_{(j)}$ consists
partially of an operation $Fix$ which has  the effect of doing $Cat_{(j-1)}$ on
each of the simplicial components,  we find that
$$
Cat _{(n)}\circ c^k(A) \subset Cat _{(n+k)}(c^k(A)),
$$
and in fact the morphism
$A\rightarrow Cat_{(n)}\circ c^k(A)$ is obtained as a successive
pushout by some of the trivial cofibrations $\Sigma \rightarrow h$  used in
\cite{nCAT} to construct $Cat_{(n+k)}$. It follows that
$c^k(A)\rightarrow Cat _{(n)}\circ c^k(A)$ is a weak equivalence of
$n+k$-precats, which
implies that
$$
Cat_{(n)}\circ c^k(A)\rightarrow Cat _{(n+k)}(c^k(A))
$$
is an equivalence of $n+k$-categories. Thus $Cat_{(n+k)}(c^k(A))$ is a
$k-1$-connected $n+k$-category.

To sum up, given a monoidal object $(A,m,e)$ in $nPC$ we obtain a $0$-connected
$n+1$-category
$$
Cat_{(n+1)}(c(A))\cong Cat _{(n)}\circ c(A).
$$
If $(A,m,e)$ is symmetric
monoidal then we obtain in fact a $k-1$-connected $n+k$-category
$$
Cat_{(n+k)}(c^k(A))\cong Cat _{(n)}\circ c^k(A)
$$
for any $k\geq 2$ (and in particular we may assume $k$ large enough so as to
land in the stable range given by Theorem \ref{thmintro}).

\newparag{easier}
If $(A,m,e)$ is a monoidal (resp. symmetric monoidal) $n$-category (i.e. a
monoidal or symmetric monoidal object in the subcategory $n-Cat\subset nPC$)
then  in the above constructions $c^k(A)$ is directly an $n+k$-category and
there
is no need for the operations $Cat$.

\newparag{ineffect}
These constructions were already used, in effect, in the construction of
the $n+1$-category $nCAT$ in \cite{nCAT}. Here one takes $nCAT$ as arising from
a strict category enriched over $nPC$, whose objects are the fibrant
$n$-categories and whose morphism-objects are the internal $\underline{Hom}$
i.e. the morphism object in $nPC$ between two objects $U,V$ is by definition
$\underline{Hom}(U,V)$. This is just a multi-object version of the above
discussion for monoidal objects in $nPC$. We recover the ``monoidal'' situation
by looking at a single object: if $B$ is a fibrant $n$-category then the full
sub-$n+1$-category of $nCAT$ consisting of the object $B$ is (by the
construction
of $nCAT$) identically the same as $c(\underline{Hom}(B,B))$ where
$\underline{Hom}(B,B)$ is given its  monoidal structure by composition.

\subnumero{The Whitehead operation}

The Whitehead operation defined in
\ref{prooflemma2}, \ref{nameit}, \ref{prooflemma3} (see also
\ref{prooflemma3bis}, \ref{wh}) is a useful way of constructing $k$-uply
monoidal $n$-categories. Indeed, one can start with {\em any} $n+k$-category
$A$, choose an object $a\in A_0$, and taking ${\bf Wh}_{>k-1}(A,a)$ yields a
$k$-uply monoidal $n$-category.  Heuristically, it is the $n+k$-subcategory of
$A$ containing as $i$-morphisms only the higher identities $1^i_a$ for all
$i\leq k-1$ (but containing the full $n$-category of endomorphisms of
$1^{k-1}_a$ from $A$). This construction intervenes in Baez-Dolan
\cite{BaezDolan} \cite{cat}, for example in the construction of ``generalized
center'' (see below).

If $A$ is an $n+k$-groupoid, corresponding to an $n+k$-truncated homotopy type,
then ${\bf Wh}_{>k-1}(A,a)$ is the $k-1$-connected $n+k$-groupoid
corresponding to the $k$-th stage in the {\em Whitehead tower} (see
\cite{Whitehead} or the exposition in \cite{BottTu}) i.e. it is the
homotopy-fiber of the morphism $$
A\rightarrow \tau _{\leq k}(A).
$$
In particular, the {\em homotopy groups} (which can be defined directly in
terms of the structure of $n+k$-groupoid cf \cite{Tamsamani}) of
${\bf Wh}_{>k-1}(A,a)$ are trivial in degrees $i\leq k$, and are the same as
the homotopy groups of $A$ (based at $a$) in degrees $i>k$. This explains the
notation ${\bf Wh}_{>k}$.

\subnumero{The ``generalized center''}

Baez and Dolan suggest in \cite{cat} a notion of ``generalized center'' ${\bf
Z}(B)$  of a $k$-uply monoidal $n$-category $B$. Technically, let us take $B$ as
being the corresponding $k-1$-connected $n+k$-category, which we furthermore
assume is fibrant. Then, in terms of our terminology from above, the definition
of Baez and Dolan \cite{cat} reads
$$
{\bf Z}(B):= {\bf Wh} _{>k+1}(c(\underline{Hom}(B,B)), x).
$$
To relate this formula to the definition of ``generalized center'' given in
\cite{cat} \S 2, recall that $c(\underline{Hom}(B,B)$ is the full
sub-$n+k+1$-category of $(n+k)CAT$ consisting of only the object $B$. On the
other hand, the Whitehead operation corresponds exactly to what is said in
\cite{cat} \S 2.

\subnumero{Suspension}

We can apply the construction $\Upsilon$ (1.3) to define the ``suspension'' of
a pointed $n$-category (again, Baez and Dolan refer to such a construction
\cite{cat}).

\newparag{suspend1}
Suppose $A$ is an $n$-precat and $a\in A_0$ is an
object. Recall that $\Upsilon (\{ a\} )= I$ is the $1$-category with objects
$0,1$ and a morphism $0\rightarrow 1$. We can use coproduct with the projection
$I\rightarrow \ast$ to define the {\em suspension}:
$$
\Sigma (A,a):=
 \Upsilon (A) \cup ^{\Upsilon (\{ a\})} \ast .
$$
Note that there is an inclusion $I\subset
\overline{I}$ (where $\overline{I}$ is the $1$-category with an isomorphism
between $0$ and $1$) which is homotopic to the projection $I\rightarrow \ast$,
so by left properness \ref{proper} we get an equivalence
$$
\Sigma (A,a)
\stackrel{\cong}{\rightarrow}
\Upsilon (A) \cup ^{\Upsilon (\{ a\})}\overline{I} .
$$

\newparag{suspend2}
The suspension is an $n+1$-precat. To get an $n+1$-category apply the operation
$Cat$ \ref{opcat}. The construction $\Sigma (A,a)$ is invariant under
weak equivalences in the variable $A$ (again by left properness \ref{proper}).

\newparag{suspend3}
Let $0$ denote the unique object of $\Sigma (A,a)$.
Note of course that we have a morphism
$$
A\rightarrow [ Cat \Sigma (A,a) ]_{1/}(0,0).
$$

\newparag{suspend4}
If $A$ is an $n$-groupoid, then the suspension $Cat \Sigma (A,a)$ is {\em
$1$-groupic} (i.e. the $i$-morphisms are invertible up to equivalence for $i\geq
1$). However, the suspension may or may not be an
$n+1$-groupoid (i.e. $0$-groupic) in this case. We have the following criterion:
if $A$ is an $n$-groupoid then the suspension $Cat\Sigma (A,a)$ is an
$n+1$-groupoid if and only if $A$ is $0$-connected (i.e. the set of
equivalence classes of objects $\tau _{\leq 0}A$ has only one element).
Let $Gr$ denote the group-completion to an $n+1$-groupoid (if it isn't
already one). Then $Gr(Cat\Sigma (A,a))$ is the $n+1$-groupoid corresponding to
the topological suspension of the space which corresponds to the $n$-groupoid
$A$.

\bigskip

We return now to the consideration of an arbitrary $A$.

\begin{lemma}
\label{mdincreases}
If $(A,a)$ is a pointed $n$-category of minimal dimension $k$ then its
suspension
$(Cat \Sigma (	A,a), 0)$ is a pointed $n+1$-category of minimal dimension $k+1$.
\end{lemma}
{\em Proof:}
The minimal dimension of the map
$$
\Upsilon (\{ a\} )\rightarrow \Upsilon (A)
$$
is at least $m(A,a) +1$, cf \ref{mdupsilon}. Preservation of minimal dimension
under coproducts implies that the minimal dimension of
$$
\ast \rightarrow \Upsilon (A) \cup ^{\Upsilon (\{ a\})} \ast
= \Sigma (A,a) $$
is the same as that of the previous map, thus it is $\geq m(A,a)+1$.
\eop

The ``delooping'' operation of Theorem \ref{thmintro} may be seen as a
suspension:

\begin{lemma}
\label{itsasuspension}
Suppose $A$ is a pointed $n$-category with minimal dimension $k$ such that
$2k\geq n+2$. Then the $n+1$-precat $X$ constructed in \ref{prethm1}
is equal to $\Sigma (A,a)$, therefore the ``delooping'' $Y$ of \ref{thm1}
is equal to $Cat\Sigma (A,a)$.
\end{lemma}
{\em Proof:}
Analyzing closely the construction $\Upsilon$ we see that
$\Upsilon (A)_{p/}$ is a disjoint union of one copy of $A$ for each
$i=1,\ldots , p$ (these are indexed by the sequences of objects
$\epsilon _0,\ldots , \epsilon _p$ where $\epsilon _i=0,1$ and $\epsilon _i\leq
\epsilon _{i+1}$). Similarly, $\Upsilon (\{ a\} )$ is a disjoint union of
one copy
of $\{ a\}$ for each $i=1,\ldots , p$. Taking the coproduct by the map
$$
\Upsilon (\{ a\} )\rightarrow \ast
$$
amounts exactly to forming the coproduct used in the construction of $X$.
\eop

We have the following variant of \ref{prethm1}, \ref{thm1}, Theorem 1,
concerning the case where the minimal dimension $k$ is not necessarily big.
Basically it says that the part of $A$ which is in the stable range is preserved
by suspension.

\begin{proposition}
\label{thm1bis}
Suppose $(A,a)$ is a pointed $n$-precat of minimal dimension $m(A,a)$. Then the
morphism
$$
A\rightarrow [Cat \Sigma (A,a)]_{1/}(0,0)
$$
has minimal dimension at least $2m(A,a)$. In particular by
\ref{relnwithtruncation} this morphism induces an equivalence on truncations
$$
\tau _{\leq 2m(A,a)-2}(A)\cong
\tau _{\leq 2m(A,a)-2} [Cat \Sigma (A,a)]_{1/}(0,0).
$$
\end{proposition}
{\em Proof:}
For any $n+1$-precat $X$, view it as a simplicial object in the
category of $n$-precats. Using minimal dimension, we can say (as in
\cite{effective}) that $X$ is
{\em $(m,k)$-arranged} if the Segal map
$$
X_{m/} \rightarrow X_{1/}\times _{X_0} \ldots \times _{X_0} X_{1/}
$$
has minimal dimension at least $k+1$.  Now, the
same procedure as in \cite{effective} may be used as a refined version of the
operation $Cat$ (the reader is referred to \cite{effective} for the
explanation).
From this we obtain that if $X$ is $(m,k)$-arranged for all $m+k\leq q$ then the
morphism $X_{1/} \rightarrow Cat(X)_{1/}$ has minimal dimension at least $q-1$.
\footnote{
[{\em Erratum:} In the last paragraph of the proof of Theorem 2.1 of
\cite{effective} the morphism in question induces an isomorphism on $\pi _i$ for
$i<n-p-1$ rather than $i<n-p$ as was stated there. This comes directly from the
conventions for marking things with an $\times$ used in that argument. As a
consequence, Theorem 2.1 of \cite{effective} should have stated that one gets
an isomorphism on $\pi _i$ for $i+m<n-1$ or more precisely that the map in
question has minimal dimension $n-m$. In particular the value of $q-1$ that
we give at this point in the present text is the correct one, rather than
$q$ as one would infer from looking at the text of \cite{effective}.]
}
As a corollary, if $X$ is $(m,k)$-arranged for all
$k\leq p$ then the above applies with $q=p+2$, so the morphism $X_{1/}
\rightarrow
Cat(X)_{1/}$ has minimal dimension at least $p+1$.

Apply the above to $X=\Sigma (A,a)$ which is the same $n+1$-precat as used in
\ref{prethm1}.  Our main estimate \ref{mainth} (cf Corollary
\ref{estimate2}) implies that the Segal maps of $X$ have minimal dimension at
least $2m(A,a)$. Thus $X$ is $(m,k)$-arranged for all $k\leq 2m(A,a)-1$ so the
morphism $A=X_{1/} \rightarrow Cat(X)_{1/}$ has minimal dimension at least
$2m(A,a)$. \eop

If $2m(A,a)\geq n+2$ we recover the statement of \ref{prethm1}.

\subnumero{The free $k$-uply monoidal $n-k$-category on one generator}
We make a few remarks about another one of Baez-Dolan's conjectures
\cite{cat}. There is a unique morphism $\partial F^k \rightarrow \ast$ and we
define the coproduct of $F^k$ using this morphism, which is an $n$-precat
$$
\sigma ^k:= F^k \cup ^{\partial F^k}\ast .
$$
We call $Cat(\sigma ^k )$ the {\em free $k$-uply monoidal $n-k$-category on one
generator}.
\end{parag}

\begin{parag}
\label{universalprop}
It is not hard to see that it satisfies the requisite universal
property: a morphism $\sigma ^k\rightarrow A$ is the same thing as
specification of an object $a\in A$ and a $k$-endomorphism of the $k-1$-fold
identity map of $a$. Let $s$ be the base object of $\sigma ^k$.
\end{parag}

\begin{parag}
\label{itsanothersuspension}
We note that
$$
\sigma ^{k+1}= \Sigma (\sigma ^{k}, s),
$$
so that $\sigma ^k$ is an iteration of the suspension operation starting with
$\sigma ^0 =2\ast $. To see this, note that
\begin{eqnarray*}
\Sigma (\sigma ^k, s)&= &\Upsilon (F^k\cup ^{\partial F^k}\ast ) \cup ^{\Upsilon
(\ast )}\ast \\
&= &\Upsilon (F^k) \cup ^{\Upsilon (\partial F^k)} \Upsilon (\ast ) \cup
^{\Upsilon (\ast )} \ast  \\
&= &\Upsilon (F^k) \cup ^{\Upsilon (\partial F^k)} \ast \\
& = & F^{k+1} \cup ^{\partial F^{k+1}}\ast \\
& = & \sigma ^{k+1} .
\end{eqnarray*}
\end{parag}

\begin{parag}
\label{conjecture}
Baez and Dolan make the
following conjecture \cite{cat}: that the $n-k$-category of endomorphisms of
the $k-1$-fold identity of $s$ in $\sigma ^k$, is the Poincar\'e $n-k$-category
$\Pi _{n-k}(X_k)$ of the space
$$
X_k = \coprod _{\ell = 0}^{\infty} C(k)_{\ell} /S_{\ell}
$$
where $C(k)_{\ell}/S_{\ell}$ is the configuration space of $\ell$ distinct
unordered points in $\rr ^k$  (\cite{cat} \S 4).
\end{parag}

\begin{parag}
\label{sketch}
Baez and Dolan already give a
sketch of an argument for this conjecture in  \cite{cat}. They point out that
for an operad $O = \{ O_{\ell }\}$, May  \cite{May} constructs the free
$O$-algebra on one point as
$$
\coprod _{\ell = 0}^{\infty} O_{\ell}/S_{\ell}.
$$
Applied to the ``little $k$-cubes'' operad $C(k)$, this gives the space $X_k$
defined above. In \cite{cat} it is argued that since $C(k)$-algebras
are $E_k$-spaces i.e. spaces with $k$-fold delooping, the free $C(k)$-algebra on
one point should be the same as the free $k$-uply monoidal $\infty$-groupoid.

This correspondence may be made precise using the results of Dunn \cite{Dunn},
which takes us closer to a rigorous proof of the conjecture \ref{conjecture}.
In effect, Dunn compares different $k$-fold delooping machines; and his model
for the $k$-fold version of Segal's machine is the same as a Segal $k$-category
(see \cite{descente} for the definition) with only one object in degrees $<k$.
Applying the $\Pi _n$ construction \cite{Tamsamani} in the top simplicial degree
(again see \cite{descente} for more details on this operation) we obtain a
correspondence between $n+k$-categories with one object in degrees $<k$ and
which
are $k$-groupic (i.e. the $i$-arrows are invertible up to equivalence for
$i>k$),
and $n$-truncated $E_k$-spaces for Dunn's Segal-machine. In
particular the $n$-truncation of the free $E_k$-space on one point  (for Dunn's
Segal-machine) is the $n+k$-category $\sigma ^k$ defined above. Now the only
thing we need to know  is that in Dunn's  comparison between different machines
for $E_k$-spaces \cite{Dunn}, the free $E_k$-spaces on one point are the same.
This should
follow directly from \cite{Dunn} using the universal properties
of the free objects, but I haven't made precise the argument. Assuming this, we
would get that the $n$-truncations of the free Segal $E_k$-space and the free
$C(k)$-algebra are the same, thus that $\sigma ^k$ is the $k$-fold delooping of
$\Pi _n(X_k)$.
\end{parag}

\begin{parag}
\label{bfsv}
The recent preprint of
Balteanu, Fiedorowicz, Schwaenzl, and Vogt \cite{BFSV}
proves a result similar to \ref{conjecture} albeit with a somewhat different
definition to start with.

\bigskip

\pagebreak[4]

\subnumero{Cohomological twisting}

\newparag{cotwist1}
In this section, let $\Xx$ be a Grothendieck site. We use the theory of
$n$-stacks on $\Xx$, see \cite{descente}. An $n$-stack $A$ is {\em
$k-1$-connected} if the truncation $\tau _{\leq k-1}(A)$ is weakly equivalent to
$\ast$. Recall that the prestack truncation will not in general be a stack, and
what we require here is that the stack associated to the prestack truncation
should be trivial. Because of this phenomenon, it is quite possible for a
$k-1$-connected $n$-stack to have global sections $\Gamma (\Xx , A)$ which are
not even $0$-connected.

\newparag{cotwist2}
Let $A^{int,0}\subset A$ be the ``interior'' of $A$ \cite{descente}, where we
retain only the $i$-morphisms which are invertible up to
equivalence. This is an $n$-stack of $n$-groupoids. Note that $\pi _0\Gamma
(\Xx ,
A^{int, 0})= \pi _0\Gamma (\Xx , A)$, so we can interpret the elements of $\pi
_0\Gamma (\Xx , A)$ as the nonabelian cohomology classes of $\Xx$ with
coefficients in the $n$-stack of groupoids $A^{int, 0}$. This leads
to the idea of ``cohomological twisting''.

\newparag{cotwist3}
We note first of all that a presheaf of $k-1$-connected $n+k$-categories on
$\Xx$ has an associated $n+k$-stack which is again $k-1$-connected. Thus, one
good way of obtaining a starting point is just to take any presheaf of
$k-1$-connected $n+k$-categories. (Or actually, the construction we do below
will also make sense if you start with {\em any}  presheaf of $n+k$-categories).
For example one could start with the constant presheaf $\underline{U}$ whose
values are a fixed $k-1$-connected $n+k$-stack.

\newparag{cotwist4}
Let $A$ be the associated $n+k$-stack. Choose an object $\alpha \in \Gamma
(\Xx ,
A)_0$, and look at
$$
{\bf Wh}_{>k-1}(\Gamma (\Xx , A),\alpha ).
$$
This is a $k-1$-connected $n+k$-category.
The idea to get new examples is to take a different base-object $\alpha \in
\Gamma
(\Xx , A)$ corresponding to a nontrivial nonabelian cohomology class
of $\Xx$ with coefficients in the interior $A^{int, 0}$. We call the resulting
$k-1$-connected $n+k$-category the {\em cohomological twist} of $A$ with respect
to the class $\alpha$.

\newparag{cotwist5}
To make things somewhat more concrete,  start  with a
$k-1$-connected $n+k$-category $U$, and let $A=\underline{U}$ be the $n+k$-stack
associated to the constant prestack with values $U$.  The interior $U^{int,
0}$ is an $n+k$-groupoid corresponding to an $n+k$-truncated space, and has a
Postnikov tower. The top stage in this tower is the morphism
$$
I:={\bf Wh}_{>n+k-1}(U^{int , 0}, 0) \rightarrow U^{int, 0},
$$
and
$$
I= K(\pi , n+k)
$$
for some group $\pi$ (abelian if $n+k\geq 2$). Thus
$$
\pi _0\Gamma (\Xx , \underline{I}) = H^{n+k}(\Xx , \pi )
$$
and this maps to $\pi _0\Gamma (\Xx , \underline{U})$. Thus we can ``twist'' by
a class in $H^{n+k}(\Xx , \pi )$. The group $\pi$ may be characterized as the
group of invertible elements in the monoid of $n+k$-endomorphisms of
$1_0^{n+k-1}$, which could be denoted $Hom (1_0^{n+k-1},1_0^{n+k-1})$
(in other words $\pi $ is the group of $n+k$-automorphisms of $1_0^{n+k-1}$).

\newparag{cotwist6}
The discussion of \ref{cotwist5} generalizes in an obvious way to the case
where we consider a presheaf of $n+k$-categories $A$ over $\Xx$. In this case
$\pi$ is a sheaf of groups on $\Xx$ and we can ``twist'' by a class in the
cohomology with coefficients in this sheaf.

\newparag{cotwist7}
We can describe what happens in a first case. Suppose $G$ is a sheaf of
groups on $\Xx$ (which we can think of as a $1$-monoidal $0$-stack), and let
$A=BG=K(G,1)$ be the $0$-connected $1$-stack associated to it. Cohomological
twisting means that we take an object $\eta \in \Gamma (\Xx , BG)$ (i.e. a
$G$-torsor or ``principal $G$-bundle'' over $\Xx$) and  look at the full
subcategory of $\Gamma (\Xx , BG)$ containing only one object $\eta$. This
subcategory is of the form $K(Ad(\eta ), 1)$, in other words it is the
$0$-connected $1$-category corresponding to the group $Ad(\eta )$ of
automorphisms of the principal $G$-bundle $\eta$. Thus, to sum up in this case,
the cohomology classes are the principal $G$-bundles, and the ``cohomological
twist'' of $G$ by a principal bundle $P$ is the group $Ad(P)$. This
classical example is important in gauge theory.

I haven't investigated what happens in any other concrete examples, but that
would be an interesting question to look at.

\end{parag}

\end{document}